\documentclass[3p]{elsarticle}

\usepackage{enumitem}
\usepackage{nameref}
\usepackage{indentfirst}
\usepackage{url}
\usepackage{graphicx}
\usepackage{listings}
\usepackage{amsmath}

\usepackage[english]{babel}
\usepackage[autostyle, english = american]{csquotes}
\MakeOuterQuote{"}
\usepackage{blindtext}
\usepackage{tablefootnote}
\usepackage{float} 
\usepackage{hyperref} 

\newcommand{\carre}{\hfill $\diamondsuit$}

\newtheorem{lemma}{Lemma}

\def\namedlabel#1#2{\begingroup
    #2%
    \def\@currentlabel{#2}%
    \phantomsection\label{#1}\endgroup
}

\begin{document}

\title{Autocorrelation Function Characterization of Continuous Time Markov Chains}

\author[mec]{G.~Rama Murthy }
\ead{rama.murthy@mechyd.ac.in}

\author[mcm]{D.G.~Down \corref{cor1}}
\ead{downd@mcmaster.ca}

\cortext[cor1]{Corresponding author}

\address[mec]{Department of Computer Science, Mahindra Ecole Centrale (MEC), Hyderabad, India}

\address[mcm]{Department of Computing and Software, McMaster University, Hamilton, Canada}

\begin{abstract} We study certain properties of the function space of
autocorrelation functions of Unit Continuous Time Markov Chains (CTMCs). It is shown
that under particular conditions, the $L^p$ norm of the autocorrelation function of
arbitrary finite state space CTMCs is infinite. Several interesting inferences are made for
point processes associated with CTMCs/ Discrete Time Markov Chains (DTMCs).
\end{abstract}

\maketitle

\section{Introduction}

Many natural/artificial phenomena are endowed with non-deterministic dynamic behavior. Stochastic/random processes with finite/countable/uncountable state space are utilized to model such dynamic phenomena. Random processes which are stationary and whose state space consists of two ($\{ +1, -1 \}$) values arise naturally in many applications e.g.\ bit transmission in communications systems. The characterization of  the autocorrelation function of such unit stochastic processes is considered an important problem \cite{Mcm}. Several interesting properties of such autocorrelation functions are studied in \cite{She,Mas}.

Wide sense stationary (or even strictly stationary) random processes naturally arise as stochastic models in a variety of applications. They also arise in time series models (AR, ARMA processes) of natural/artificial phenomena. In most such models, the autocorrelation function, $R(\tau)$ is integrable and hence the power spectral density (the Fourier transform of $R(\tau)$) exists. There are a number of applications (such as random telegraph signals) where the underlying stochastic process can be in only one of two states, so-called unit random processes. An example can be taken from detection theory. Suppose that a random process, $x(t)$ is provided as input to a threshold detector, where the output $y(t)$ is the sign of $x(t)$, i.e., $y(t)=+1$ if $x(t) \geq 0$ and $y(t)=-1$ if $x(t)<0$. Thus the output is a unit process and can be shown to be Markovian under some conditions on $x(t)$. In this case, the autocorrelation of $y(t)$ is in general not integrable.  More generally, quantization of a general random process leads to a finite-state random process, which is often Markovian, so the study of more general finite-state processes is also of interest. With this in mind, we are motivated to study the function space of finite state Markovian processes. Masry \cite{Mas} has studied the functional space properties of stationary unit random processes. However, the study of the integrability of $R(\tau)$ was not undertaken.  To the best of our knowledge, integrability, or more generally the $L^p$-norm of $R(\tau)$ for finite state Continuous Time Markov Chains (CTMCs) has not ben investigated. In this paper, we determine conditions under which the autocorrelation function is not integrable, and by extension conditions under which the $L^p$-norm of $R(\tau)$ approaches zero as $p\rightarrow \infty$. 

To put our work into context, there is related work in three directions: the characterization of autocorrelation functions of random processes, the characterization of point processes, and the use of autocorrelation properties in the analysis of stochastic models, in particular the analysis of queues. For the characterization of autocorrelation functions, we point the reader to work in time series analysis \cite{deb11,deglac03,ino08} and in telecommunications \cite{eunjinhonson04}. Point process characterization has been studied in \cite{hagvanmol99,hua90,ivamerpla07}. Properties of autocorrelation functions have been employed to determine appropriate simulation strategies for queues \cite{whi91} and is a feature of modelling arrival traffic to queues, using Markovian Arrival Processes, see \cite{klelinloh03}, for example.

This paper is organized as follows. In Section~\ref{s2}, the autocorrelation function of a Unit ($\{ +1, -1 \}$ state) CTMC is computed and the structure of the function space is studied. In Section~\ref{s3}, the autocorrelation function of a symmetric state space CTMC (i.e., the state space is $\{ -N, -N+1,\ldots,-1, +1, \ldots+N-1, +N \}$) is computed and the finiteness of its $L^{p}$ norm is discussed. It is shown that under some conditions, the autocorrelation function is not integrable. In Section~\ref{s4}, various interesting inferences are made for point processes. Finally, the paper concludes in Section~\ref{s5}.
  
  \section{Auto-Correlation Function of Homogeneous Unit CTMC: Integrability}
  \label{s2}

In the following, we consider a homogeneous Continuous Time Markov Chain (CTMC) (with the state space being $\{ +1, -1 \}$) whose generator matrix is denoted by $Q$. We assume that the resulting stochastic process is wide sense stationary (and not necessarily strictly stationary).

The autocorrelation function is given by
\[ R(t,t+\tau) = E [X(t)X(t +\tau)]. \]
Since the CTMC is homogeneous, we have that
\[ R( t,t+\tau) = R(\tau) = E[X(0)X(\tau)]. \]
Since the CTMC is a unit random process (i.e., the state space is $\{ +1,-1 \}$), we have that
\[ R(\tau) = P \{ X(0) = X(\tau) \} - P \{ X(0)  \neq X(\tau) \}. \]
Combining this with
\[ P \{ X(0)  \neq X(\tau) \} = 1 - P\{ X(0) = X(\tau) \}  \]
yields
\begin{equation}
\label{e1}
R(\tau) = 2P\{ X(0) = X(\tau) \}  -1.
\end{equation}
It remains to compute $P \{ X(0)=X(\tau) \}$.
\begin{eqnarray}
\lefteqn{P\{ X(\tau) = X(0) \}} \nonumber\\
& = & P\{ X(\tau) = +1,X(0) = +1\} + P\{ X(\tau) = -1,X(0) = -1 \} \nonumber \\
& = & P\{ X(\tau) = +1 | X(0) = +1 \} P\{ X(0) = +1 \} + P\{ X(\tau) = -1 | X(0) = -1\} P\{ X(0) = -1 \}. \label{e2}
\end{eqnarray}
The conditional probabilities in the above expression are computed using
the transient probability distribution for a homogeneous CTMC. This computation requires the determination of the matrix exponential of $Q$, i.e., $e^{Qt}$. Using the
spectral representation, $e^{Qt}$ is computed below.

\begin{itemize}
\item Computation of $e^{Qt}$. Let the generator matrix of $\{ +1, -1 \}$ state CTMC (i.e., Unit CTMC) be given by
\[ Q = \begin{bmatrix}
   - \alpha & \alpha \\
  \beta & - \beta \\
  \end{bmatrix}.
\]
Since $Q$ is a rank one matrix, the eigenvalues are $\{ -( \alpha + \beta),0 \}$. Also, let the right eigenvectors be given by $\{ \bar{g}_1, \bar{g}_2 \}$ (i.e., column vectors) and let the left eigenvectors be given by $\{ \bar{f}_1,\bar{f}_2 \}$ (i.e., row vectors). Thus, we have that
\[ e^{Qt} = e^{-(\alpha + \beta)t} \bar{g}_1 \bar{f}_1 + \bar{g}_2 \bar{f}_2. \]
\end{itemize}

Now, we discuss the computation of eigenvectors.

\begin{itemize}
\item Computation of $\bar{g}_1$. By definition of right eigenvector, we have that
\[ Q\bar{g}_1 = \  - (\alpha + \beta) \bar{g}_1. \]
By letting $\bar{g}_1 = [ x \ y ]^T$ and using the expression for the generator matrix $Q$,
\[ \begin{bmatrix}
 - \alpha & \alpha \\
\beta & - \beta \\
\end{bmatrix} \begin{bmatrix}
x \\
y \\
\end{bmatrix} = - (\alpha + \beta)\ \begin{bmatrix}
x \\
y \\
\end{bmatrix}. \]
Solving the above linear equations, we have $y = - \frac{\beta}{\alpha} x$. Thus letting $x = 1$, we have $y = - \beta / \alpha$. Hence
\[ \bar{g}_1 = \begin{bmatrix}
1 \\
\frac{-\beta}{\alpha} \\
\end{bmatrix}. \]

\item Computation of $\bar{g}_2$. Since $Q$ is the generator of the CTMC, we have
\[ \begin{bmatrix}
   - \alpha & \alpha \\
  \beta & - \beta \\
  \end{bmatrix} \begin{bmatrix}
  1 \\
  1 \\
  \end{bmatrix}  = \begin{bmatrix}
  0 \\
  0 \\
  \end{bmatrix}. \]
Thus, the right eigenvector corresponding to the zero eigenvalue is given by
\[ \bar{g}_2 = \begin{bmatrix}
1 \\
1 \\
\end{bmatrix}. \]

\item Computation of $\bar{f}_1$. Since the generator matrix $Q$ is diagonalizable, we have
\[ \bar{f}_1 \bar{g}_1 = 1 \ \mbox{and} \ \bar{f}_1 \bar{g}_2 = 0. \]
Letting $\bar{f}_1 = [a \ b]^T$, $\bar{f}_1 \bar{g}_2  = 0$ implies $a = -b$.  Furthermore, $\bar{f}_1 \bar{g}_1 = 1$ implies that $b = -\frac{\alpha}{\alpha + \beta}$. Hence
\[
\bar{f}_1 = \left[ \begin{array}{cc}
\frac{\alpha}{\alpha + \beta} & - \frac{\alpha}{\alpha + \beta}
\end{array} \right].
\]

\item Computation of $\bar{f}_2$. Letting $\bar{f}_2 = [ c \ d]^T$, $\bar{f}_2 \bar{g}_1 = 0$ implies that $c = d \frac{\beta}{\alpha}$. Furthermore, $\bar{f}_2 \bar{g}_2 = 1$ implies $d = \frac{\alpha}{\alpha + \beta}$. Hence, we have that
\[ \bar{f}_2 = \left[ \begin{array}{cc}
\frac{\beta}{\alpha + \beta} & \frac{\alpha}{\alpha + \beta}
\end{array}  \right]. \]
\end{itemize}
Thus, we have effectively computed the matrix exponential $e^{Qt}$.

The transient probability distribution of our homogeneous, Unit CTMC is given by
$\bar{\pi}(\tau) = [P\{ X(\tau) = +1\} \ P\{ X(\tau) = -1 \}]$.
Using the matrix exponential, the transient probability distribution is
given by
\begin{eqnarray*}
\bar{\pi}(\tau) & = & \bar{\pi}(0)\ e^{Q\tau}\\
& = & \bar{\pi}(0)\ \left[ e^{-( \alpha + \beta)\tau} \bar{g}_1 \bar{f}_1 + \bar{g}_2 \bar{f}_2 \right]\\
& = & \bar{\pi}(0)\ \left[ \begin{matrix}
P\{ X(\tau) = +1 | X(0) = +1\} & P\{ X(\tau) = -1 | X(0) = +1\}  \\
P\{ X(\tau) = +1 | X(0) = -1\} & P\{ X(\tau) = -1 |X(0) = -1\} \\
\end{matrix} \right].
\end{eqnarray*}
We now provide an explicit expression for the matrix exponential $e^{Q \tau}$.
\begin{eqnarray*}
e^{Q\tau} & = & e^{-( \alpha + \beta)\tau} \begin{bmatrix}
1 \\
\frac{- \beta}{\alpha} \\
\end{bmatrix} \left[ \begin{matrix}
\frac{\alpha}{\alpha + \beta} & -\frac{\alpha}{\alpha + \beta} \\
\end{matrix}  \right] + \begin{bmatrix}
1 \\
1 \\
\end{bmatrix}\left[ \begin{matrix}
\frac{\beta}{\alpha + \beta} & \frac{\alpha}{\alpha + \beta} \\
\end{matrix}  \right]\\
& = & \begin{bmatrix}
e^{-( \alpha + \beta)\tau}\frac{\alpha}{\alpha + \beta} + \frac{\beta}{\alpha + \beta} & -e^{-( \alpha + \beta)\tau}\frac{\alpha}{\alpha + \beta}  + \frac{\alpha}{\alpha + \beta}  \\
 -e^{-( \alpha + \beta)\tau}\frac{\beta}{\alpha + \beta} + \frac{\beta}{\alpha + \beta} & e^{-(\alpha + \beta)\tau}\frac{\beta}{\alpha + \beta} + \frac{\alpha}{\alpha + \beta} \\
\end{bmatrix}.
\end{eqnarray*}
We utilize the expression for $e^{Q\tau}$ in (\ref{e2}) to compute $P\{ X(\tau) = X(0) \}$.
\begin{eqnarray}
\lefteqn{P\{ X(\tau) = X(0) \}} \nonumber\\
& = & \ P\{ X(\tau) = +1 | X(0) = +1\}  P\{ X(0) = +1 \} + P\{ X(\tau) = -1 | X(0) = - 1\}  P\{ X(0) = -1 \} \nonumber\\
& = &
\left[ e^{- (\alpha + \beta)\tau}\frac{\alpha}{\alpha + \beta} \right] P\{ X(0) = +1 \}  + \left[ e^{-(\alpha + \beta)\tau}\frac{\beta}{\alpha + \beta} \right] P\{ X(0) = -1  \} \nonumber \\
& & + \frac{\beta}{\alpha + \beta}P\{ X(0) = +1\  \}  + \frac{\alpha}{\alpha + \beta} P\{ X(0) = -1 \} \label{e3}
\end{eqnarray}

Now, we compute the equilibrium probability distribution of the unit CTMC. Let the equilibrium probability distribution be denoted by the row vector
$\bar{\pi}$. Thus
$\bar{\pi} = [\pi_{1} \ \pi_{- 1}]$ and $\overline{\pi}Q = 0$. Solving
the linear equations, we have $\pi_{1} = \frac{\beta}{\alpha + \beta}$ and $\pi_{- 1} =
\frac{\alpha}{\alpha + \beta}$.

\textbf{Note:} We now assume that the initial probability distribution
equals the equilibrium probability distribution, i.e., $\overline{\pi}(0) = \overline{\pi} = [\pi_{1} \ \pi_{- 1}]$.
In this case, the transient probability distribution is also equal to the equilibrium probability distribution, which allows us to write
\[ P \{ X (\tau) = X(0) \} = \frac{\alpha^{2} + \beta^{2}}{(\alpha + \beta)^{2}}. \]
Now substituting for $P\{ X(\tau) = X(0) \}$
in (\ref{e1}), we have that
\[ R(\tau) = 2P\{ X(\tau) = X(0)  \}  - 1  = 2\left[ \frac{\alpha^{2} + \beta^{2}}{(\alpha + \beta)^{2}} \right] - 1. \]
Simplifying, we have that
\begin{equation} \label{e:R}
R(\tau) = \left[ \frac{(\alpha - \beta)^{2}}{(\alpha + \beta)^{2}} \right].
\end{equation}
Thus, the autocorrelation function is a constant when the initial
probability distribution equals the equilibrium probability
distribution. Further, if $\alpha = \beta$, we have that $R(\tau)=0$, i.e., it is identically zero. Since the initial probability distribution is assumed to equal the
equilibrium probability distribution, with $\alpha = \beta$, we have
\[ E[X(0)] = (1)\frac{\beta}{\alpha + \beta} +(-1) \frac{\alpha}{\alpha + \beta} = 0, \]
i.e., with mean zero initial random variable, (\ref{e:R}) holds for
the autocorrelation function.

Suppose, we consider the case where the initial probability distribution
is not equal to the equilibrium probability distribution and is
arbitrary, i.e.,
\[ P\{ X(0) = + 1\} = q \quad \mbox{and} \quad P\{ X(0) = -1\} = 1 - q. \]
Substituting in (\ref{e3}), we have that
\begin{eqnarray*}
\lefteqn{P\{ X(\tau) = X(0)\} }\\
& = & \left[ e^{-(\alpha + \beta)\tau}\frac{\alpha}{\alpha + \beta} \right] q  + \left[ e^{-(\alpha + \beta)\tau}\frac{\beta}{\alpha + \beta} \right] (1-q) + \frac{\beta}{\alpha + \beta} q  + \frac{\alpha}{\alpha + \beta} (1 - q)\\
& = & 
e^{-(\alpha + \beta)\tau}\ \left( \frac{\alpha - \beta}{\alpha + \beta} \right)q + e^{-(\alpha + \beta)\tau}\left( \frac{\beta}{\alpha + \beta} \right) + \left( \frac{\beta - \alpha}{\alpha + \beta} \right)q + \frac{\alpha}{\alpha + \beta}.
\end{eqnarray*}
Substituting in $R(\tau) = 2P\{ X(\tau) = X(0) \}  - 1$,
\[ R(\tau) = 2\left[ e^{-(\alpha + \beta)\tau} \left( \frac{\alpha - \beta}{\alpha + \beta} \right)q + e^{-(\alpha + \beta)\tau}\left( \frac{\beta}{\alpha + \beta} \right) \right] + 2 \left[ \left( \frac{\beta - \alpha}{\alpha + \beta} \right)q + \frac{\alpha}{\alpha + \beta} \right] - 1. \]
Thus
$R(\tau)$ is a sum of a function of $\tau$ and a constant term, which we write as $R(\tau)=f(\tau)+c$, where
\[ c = 2\left[ \left( \frac{\beta - \alpha}{\alpha + \beta} \right)q \right] + \frac{\alpha - \beta}{\alpha + \beta}. \]
If $q = \frac{1}{2}$, we have that $c$ is zero and thus
\[
R(\tau) = e^{-(\alpha + \beta)\tau}\ \left[ \left( \frac{\alpha - \beta}{\alpha + \beta} \right) + 2\left( \frac{\beta}{\alpha + \beta} \right) \right]
= e^{-(\alpha + \beta)\tau}.
\]
Now, if $\alpha = \beta$, and since the
autocorrelation function is symmetric, we have
\[ R(\tau) = e^{-2\alpha |\tau|}. \]


As discussed earlier, an interesting problem is to
characterize the autocorrelation function of Unit CTMCs, i.e., $\{ +1, -1 \}$-state CTMCs.
The following lemma presents one possible characterization of the function
space of autocorrelation functions of Unit CTMCs

\begin{lemma}
Consider a Unit CTMC with $\alpha \neq \beta$. The
autocorrelation function of such a Unit CTMC is not in
$L^{p}[R(\tau)]$ for any $p \geq 1$ (the $L^p$ norm of the autocorrelation function is infinite) . However, as $p$ tends to $\infty$, the $L^{p}$ - norm of the autocorrelation function, $R(\tau)$, approaches a finite constant. Further the $L^{\infty}$-norm is equal to one.
\end{lemma}

{\bf Proof.} In the case of unit CTMCs with $\alpha \neq \beta$ and the initial probability distribution is equal to the equilibrium probability distribution,
the expression for autocorrelation function is identically (for all
time) a non-zero constant which is less than one in  magnitude. Thus, such a function is not in $L^{p}(R)$ for any
$p \geq 1$. Thus, such an autocorrelation function is necessarily not integrable.
Since the constant term is less than one in magnitude, it readily
follows that the $p^{th}$ power of it approaches zero. Hence the
$L^{p}$-norm of the autocorrelation function, $R(\tau)$, approaches a finite constant.

In the case of arbitrary initial probability distribution, if
$q \neq \frac{1}{2}$, then $c$ is non-zero. Further if
$q \leq \frac{1}{2}$, $0 \leq c \leq 1$ and
$q \geq \frac{1}{2}$, $-1 \leq c \leq 0$. Thus, even in
this case, the autocorrelation function of such a Unit CTMC is not in
$L^{p}(R)$ for any $p \geq 1$. Hence,
it is necessarily not integrable. Also, as the constant term is strictly
less than one in magnitude, its $p^{th}$ power approaches zero, the $L^{p}$-norm of the autocorrelation function, $R(\tau)$, approaches a finite constant. \carre

In the following discussion, we generalize the above results. We
consider symmetric state space CTMCs (with state space $\{ -N,
-N+1,\ldots,-1, +1, \ldots, +N-1, +N \}$) as well as arbitrary state space CTMCs. It is shown that the existence
of an equilibrium probability distribution ensures that the expression for
autocorrelation function has a constant part that is not
zero under some conditions.

\section{Auto-Correlation Function of Homogeneous Finite
  State Space CTMC}
\label{s3}

\subsection{Asymmetric State Space}

We now prove that for any finite state space CTMC, the autocorrelation function is not integrable and in
fact the $L^{p}$-norm of the autocorrelation function, $R(\tau)$,
is infinite for any $p \geq 1$. Without loss of generality, the state space of the CTMC is assumed to be $\{ +1, +2, \ldots, +N \}$.

For a homogeneous CTMC, we have that $R(t,t + \tau) = R(\tau) = E[X(0)X(\tau)]$. Since the CTMC under consideration has finite state space, we have that
\begin{eqnarray*}
R(\tau) & = & \sum_{i = 1}^{N}\sum_{j = 1}^{N} ij P\{ X(0) = i,X(\tau)=j  \}\\
& = & \sum_{i = 1}^{N}\sum_{j = 1}^{N} ij P\{ X(\tau) = j | X(0) = i \} P\{ X(0) = i \}\\
& = & \sum_{i = 1}^{N}\sum_{j = 1}^{N} ij \left( e^{Q\tau} \right)_{ij} P\{ X(0) = i  \}.
\end{eqnarray*}
But, we have that
\[ e^{Q\tau} = \sum_{k = 1}^{N}e^{\gamma_{k}\tau} E_{k}, \]
where $E_{k}$ is the residue matrix such that $E_{k} = \bar{f}_k\ \bar{g}_k$ with $\bar{f}_k$ being the right eigenvector of $Q$ corresponding to the eigenvalue $\gamma_k$ and $\bar{g}_k$ being the left eigenvector of $Q$ corresponding to the eigenvalue $\gamma_k$. Let
$|\gamma_1| \geq |\gamma_{2}| \geq \cdots \geq | \gamma_{N} | = 0$ (since $Q$ is the generator matrix, we have that $\gamma_{N} = 0$).

\begin{eqnarray}
R(\tau) & = & \sum_{i = 1}^{N} \sum_{j = 1}^{N} ij \left( e^{Q\tau} \right)_{ij} P\{ X(0) = i \} \nonumber \\
& = & \sum_{i = 1}^{N} \sum_{j = 1}^{N} ij \left[ \sum_{k = 1}^{N}e^{\gamma_{k}\tau} E_{k}\ \right]_{ij} P\{ X(0) = i  \} \nonumber \\
& = & \sum_{i = 1}^{N} \sum_{j = 1}^{N} ij P\{ X(0) = i \} \left[ \sum_{k = 1}^{N - 1}e^{\gamma_{k}\tau} E_{k} \right]_{ij}
+ \sum_{i = 1}^{N}\sum_{j = 1}^{N}
ij P\{ X(0) = i  \}[ \bar{f}_N \ \bar{g}_N]_{ij} \label{e:star}
\end{eqnarray}
We would like to write (\ref{e:star}) as $f(\tau)+c$, where
\[
c = \sum_{i = 1}^{N}\sum_{j = 1}^{N} ij P\{ X(0) = i \} [\bar{f}_N \ \bar{g}_N]_{ij}.
\]
We must justify that $c$ is not a function of $\tau$. To do this, it
is first noted that $\bar{f}_N$ is the equilibrium
probability vector of the CTMC, i.e., $[ \pi_{1} \ \pi_{2} \cdots \pi_{N}]$.
It is also the left eigenvector of the generator matrix $Q$ corresponding to the
eigenvalue zero. Also
$\bar{g}_N = [ 1 \ 1\cdots 1]^{T}$, i.e.,
a column vector of all ones. It is the right
eigenvector of the generator matrix $Q$ corresponding to the eigenvalue zero. By the Perron-Frobenius Theorem, all of the components of the vectors $\{
\bar{f}_N, \bar{g}_N \}$ are non-negative. Thus,
$c$ involves a sum of non-negative quantities and does not depend on
$\tau$. It readily follows that $[ \bar{f}_N \ \bar{g}_N]_{ij} = \pi_{j}$. Let $P\{ X(0) = i \} = q_{i}$. Hence, we have that
\begin{eqnarray*}
c & = &  \sum_{i = 1}^{N}\sum_{j = 1}^{N} ij q_{i} [\bar{f}_N \ \bar{g}_N]_{ij} \\
& = & \left( \sum_{i = 1}^{N} iq_{i} \right) \left( \sum_{j = 1}^{N} j\pi_{j} \right) \\
& = & E[X(0)]E[Z],
\end{eqnarray*}
where $Z$ is the random variable associated with the equilibrium
probability distribution. Thus, $c$ is zero if and only if either $E[X(0)]$ or $E[Z]$ is zero (or both). For our choice of state space, this is not possible.

\textbf{Note:} If the initial probability distribution equals the
equilibrium probability distribution, then $c = E[X(0)]^2 = E[Z]^2$.

\textbf{Note:} In general, the state space of a finite state space CTMC
could assume negative values (as in the case of a Unit CTMC). It can be
readily verified that if the state space of CTMC is the set $\{ -N, -N+1,
\ldots, -1, +1, +2,\ldots, +N \}$ and the equilibrium probability
distribution is the uniform distribution, then $c$ equals zero.
Thus, such a condition is sufficient for $c$ to be zero. We
conjecture that the condition is also necessary for $c$ to be zero.


Now, $f(\tau)$ is a sum of decaying exponentials. It can be easily
verified that $f(\tau)$ is integrable. More
generally, $f(\tau)$ corresponds to a function which is in
$L^{p}(R)$ for $p \geq 1$. But, when $c$ is
non-zero, then
$R(\tau)$ is not in $L^{p}(R)$
for any $p \geq 1$. 
If $c \neq 0$, then $\int (R(\tau))^{p} d\tau$ is infinite for every $p \geq 1$. Further if $|c| < 1$, then $\int (R(\tau))^{p}d\tau$ approaches zero as $p \rightarrow \infty$. The following lemma readily follows from the above discussion.

\begin{lemma}
Consider an arbitrary finite state CTMC with $c$
being non-zero. The autocorrelation function of such a unit CTMC is
not in $L^{p}(R)$ (where $R$ is the
real valued lag) for any $p \geq 1$,  i.e., such an autocorrelation
function is necessarily not integrable. However, as
$p$ tends to $\infty$, the $L^{p}$-norm of the autocorrelation function,
$R(\tau)$ approaches a finite constant.
\end{lemma}

{\bf Note:} We can provide further details about $c$, in the case that the
random process is a Unit CTMC.

Case (A): $X(0) = Z$, i.e., the initial probability distribution equals the
equilibrium probability distribution.
\[
c = E[X(0)]^2 = E[Z]^2 = (1)\frac{\beta}{\alpha + \beta}\  + (-1)\frac{\alpha}{\alpha + \beta} = \left( \frac{\alpha - \beta}{\alpha + \beta} \right)^{2}.
\]

CASE (B) : $X(0) \neq Z$, i.e., the equilibrium probability distribution is arbitrary, i.e., $P(X(0) = +1) = q$ and $P(X(0) = -1) = 1 - q$. Hence, we have that
\begin{eqnarray*}
E[X(0)] & = & q(1) + (1 - q)(-1) = 2q-1\\
E[Z] & = & \frac{\beta - \alpha}{\alpha + \beta}.
\end{eqnarray*}
Hence $E[X(0)]E[Z] = 2\left[ \left( \frac{\beta - \alpha}{\alpha + \beta} \right)q \right] + \frac{\alpha - \beta}{\alpha + \beta}$
which is consistent with our earlier result.

\textbf{Claim:} We now interpret the above results for arbitrary
Continuous Time Markov Chains (CTMCs) for which an equilibrium
distribution exists. Since $R(\tau) = E[X(0)X(\tau)]$,
we have that the following asymptotic result holds true.
\[ \lim_{\tau \rightarrow \infty} R(\tau) = E[X(0)]E[Z], \]
where $Z$ is the equilibrium random variable.

This result agrees with the fact that asymptotically the initial random
variable $X(0)$ and the equilibrium random variable are independent. Also, these results can be easily generalized to countable state space CTMCs.

\section{Finite State Space Continuous Time Markov Chains: Point
  Processes}
  \label{s4}


It is well known that the interarrival times of a Poisson process are exponentially distributed random variables. Also, the sojourn times in every state of a finite state CTMC are exponentially distributed random variables. This observation has been explored in \cite{cin75}, for example, to establish that when successive visits to a state of a CTMC are stitched together, a Poisson process naturally results. Hence, an arbitrary finite state CTMC can be viewed as a superposition of point processes. From a practical viewpoint, the superposition of point processes naturally arises in applications, such as packet streams in packet multiplexers. Such packet streams have been modelled in \cite{sriwhi86}, for example. Several versatile point processes have also been studied in \cite{hefluc86,neu79}, amongst others. Such Markovian point processes are actively utilized in queueing theoretic applications. In discrete time, sojourn times in any state of a DTMC are geometrically distributed random variables. In this sense, a DTMC can be viewed as a superposition of discrete time point processes (which could be dependent). Thus, homogeneous CTMCs/DTMCs provide interesting models  of superpositions of point processes. In other words, many well known results related to equilibrium/transient probability distributions of such Markov chains, effectively provide novel inferences related to the superposition of packet streams in communication networks.
Our main goal in this section is to  to derive several such inferences, considering CTMCs to illustrate the idea.

\subsection{Equilibrium/Transient Inferences related to Point Processes}

  Let $\overline{\pi} =
  [\pi_{1} \ \pi_{2} \cdots \pi_{N} ]$ be the
  vector of equilibrium probabilities.

\begin{enumerate}
\item
  Probability that in equilibrium the arrival is from
  the $j^{th}$ Point Process = $\pi_{j}$, for $1 \leq \ j \leq N$.
\item Probability that at time $\tau$ (under transient conditions), the
  arrival is from the $j^{th}$ point process, $P\{ X(\tau) = j  \} = [\overline{\pi}(0)e^{Q\tau}]_{j}$,
 i.e., $j^{th}$ component of such a vector for $1 \leq \ j \leq N$.
\item
  Conditional Probability that at time $\tau$ (under transient
  conditions), the arrival is from the
  $j^{th}$ point process given that at time zero  the  process was  in state $i$ is equal to
$[e^{Q\tau} ]_{ij}$ for $1 \leq  i \leq N$ and $1 \leq  j \leq N$.

  Many such inferences can be made related to the point processes in
  superposition given the interesting quantities associated with the
  transient/equilibrium probability distribution of corresponding DTMC/CTMC models.
\end{enumerate}

\section{Conclusion}
\label{s5}

We have computed the autocorrelation function of a unit CTMC
and the conditions for integrability (more generally finiteness of the
$L^{p}$-norm) were established. More generally, the function space structure of
  arbitrary finite state space CTMC was explored. Interesting inferences
  related to point processes (in a superposition point process) were
  made based on their relationship to finite state space Markov chains.
  
\vspace{.1in}

\noindent {\bf Acknowledgment.} The first author is supported by internal funding for research and development from Mahindra \'{E}cole Centrale, Hyderabad (for the year 2018-19). The second author is supported by NSERC under the Discovery Grant program.


\end{document}